
 \input amstex
\documentstyle{amsppt}
\magnification=\magstep1
\define\CalC{{\Cal C}}
\define\bfR{{\bold R }}
\define\bfZ{{\bold Z}}
\define\bfC{{ \bold C}}
\define\CalB{{ \Cal B}}
\define\done{{df(x) \over dx}}
\define\dtwo{{d^2 f(x) \over dx^2}}
\define\dthree{{d^3 f(x) \over dx^3}}
\define\dfour{{d^4 f(x)\over dx^4}}
\define\dfive{{d^5 f(x)\over dx^5}}

\define\dseven{{d^7 f(x)\over dx^7}}

\define\ski{\cr &\qquad+}

\topmatter
\title The Toda Hierarchy and the KdV Hierarchy\endtitle

\author D. Gieseker \endauthor
\address Department of Mathematics, UCLA \endaddress
\email dag\@math.ucla.edu \endemail
\thanks I wish to thank Russ Caflisch and Stan Osher for some helpful
conversations. \endthanks
\thanks Partially supported by NSF Grant DMS 93-05657\endthanks
\endtopmatter
\document

McKean and Trubowitz \cite{2} showed that the theory of the KdV
equation  $${\partial \over \partial t}g(x,t)= {\partial^3\over
\partial x^3}g(x,t)-6g(x,t){\partial g\over \partial x}(x,t).$$
is intimately related to the geometry of a related  hyperelliptic curve of
infinite genus, the Bloch spectrum $\CalB_{g_{t}}$ of the operator
$$L_{g_{t}}:\psi\to {d^2\over dx^2}\psi(x)+g(x,t)\psi(x),$$ where
$g_t=g(x,t).$ As was known
classically, $\CalB_{g_{t}}$ is independent of $t,$ when $g(x,t)$ evolves
according to the KdV equation.  Our purpose in this paper is to develop
a theory of finite difference operators and their Bloch spectra and
isospectral flows which  mimics the KdV theory.   The basic idea of
this paper is  to use the theory of the periodic Toda chain of length
$N$.  Here again, the periodic Toda chain can be understood in terms of
a finite genus hyperelliptic curve and isospectral deformations, as
van Moerbeke  discovered. For instance, see \cite{3}. So one
would like to see what the relation of the Toda hierarchy is to the KdV
hierarchy, how the conserved quantities are related and so forth. A
start on these matters has been obtained by Toda in \cite{4}.  In this
paper, the idea is that if we choose the initial data for the periodic
Toda chain very carefully, then the evolution of this data under the
various equations of the Toda hierarchy looks similar to the evolution
of $f$ under the KdV
hierarchy. Given $f,$ we will find a canonical choice of the initial
data of the Toda equations so that the flow of this initial data under
the Toda hierarchy looks like the flow of $f$ under the KdV hierarchy, at least
to high accuracy.  This choice will be given by an asymptotic series
in $N^{-1}.$  The main result of this paper is the formulation and
partial verification of the Conjecture given below. The method studied here
also gives a way of producing
analytically  approximate solutions to the Toda chain hierarchy, at
least conjecturally.     My motivation is to use this case as a model
for constructing finite genus models for the KP hierarchy, which I have
studied in \cite{1}. I also hope to use these methods to develop
discrete models of the sine-Gordon equation, the non-linear Schr\"odinger
equation and other infinite
dimensional integrable systems that are related to isospectral flows. At
the end of this paper, I discuss using the methods developed in this
paper to construct a
finite difference equation to numerically solve the  KdV equation.  This
finite difference equation has many conserved quantities.

We first review the KdV hierarchy from an isospectral point of view.
Let $T$ be the translation operator defined on real  valued functions
on $\bfR$ by $$T(g)(x)=g(x+1).$$  Let $\CalC$ be the set of functions
$g$ so that $g$ is analytic on $\bfR$ and $T(g)=g,$ i.e. $g$ is
periodic with period 1.  If  $g\in \CalC$, then we can define an
operator $$L_g(\psi)=-\psi''+g(x)\psi(x)$$ on the space of infinitely
differentiable functions $\psi$ on $\bfR.$  We define the Bloch
spectrum $\CalB_g$ of $L_g$ to be the set of $(\lambda,\alpha)\in \bfC
\times \bfC^*$ so that there is a non-zero function $\psi$ with
$L_g(\psi)=\lambda \psi$ and $T(\psi)=\alpha \psi$.  The KdV hierarchy is a
sequence of non-linear differential operators $D_i:\CalC\to \CalC$ so
that for any $i$ if we have an analytic  function $g(x,t)$ and we
define $g_t(x)$ to be $g(x,t)$  and $$D_i(g_t)={d \over dt} g_t,$$ then
$\CalB_{g_t}$ is independent of $t.$

The Toda hierarchy is quite similar to the KdV hierarchy.  Let us fix
an positive integer $N.$  If $A:\bfZ \to \bfC$, then we can again define
the translation operator $T$ by $T(A)(n)=A(n+N).$ Let $\CalC_N$ be the
set of $(A,B)$ with $T(A)=A$ and $T(B)=B.$  Given $(A,B) \in \CalC_N$,
we can form an operator on the space of all functions $\psi:\bfZ\to
\bfC$ by
$$L_{(A,B)}(\psi)(n)=(-\psi(n+1)+A(n)\psi(n)+B(n)\psi(n-1))N^2.$$ We
define the Bloch spectrum $\CalB_{(A,B)}$ of $L_{(A,B)}$ to be the set
of $(\lambda,\alpha)\in \bfC \times \bfC^*$ so that there is a non-zero
function $\psi$ with $L_{(A,B)}(\psi)=\lambda \psi$ and $T(\psi)=\alpha
\psi.$  Then there is a sequence of non-linear difference operators
$D_{1,k}(A,B)$  and $D_{2,k}(A,B)$ so that if $A_t$ and $B_t$ are in
$\CalC_N$ and $A_t(n)$ and $B_t(n)$ are differentiable functions of $t$
and we have $$D_{1,k}(A_t,B_t)={d\over dt}A_t$$ and
$$D_{2,k}(A_t,B_t)={d\over dt}B_t,$$ then   $\CalB_{(A_t,B_t)}$ is
independent of $t.$

The Toda hierarchy can be defined in the following way:  We first
inductively define a  sequence of complex valued functions   $d_i$  on
$\bfZ$  as follows:  We take $d_i=0$ for $i<0$ and $d_0=1.$  We also
take $d_i(0)=0$ and finally
$$d_i(n+1)=d_i(n)+A(n)d_{i-1}(n)+B(n)d_{i-2}(n).$$ The first Toda
equation is defined by $$D_{1,1}(A,B)(n)= (B(n)-B(n+1))N$$ and
$$D_{2,1}(A,B)(n)=B(n)(A(n)-A(n-1))N.$$  To define the $k^{th}$ member
of the Toda hierarchy ($k>1$), we define a sequence of functions
$a_{k-1},\ldots,a_0$ by descending induction on $p$ by
$$a_p(n)=d_{k-p}(n+k)-d_{k-p}(n)-\sum_{r=p+1}^{k-1}
a_r(n)d_{r-p}(n+r).$$ We define
$$D_{1,k}(A,B)(n)=(a_1(n)B(n+1)-B(n)a_1(n-1))N$$ and
$$D_{2,k}(A,B)(n)=B(n)(a_0(n)-a_0(n-1))N.$$  The numbers $d_i(N)$ are
conserved quantities of the flow. If necessary, we denote  $d_i(N)$ by
$d_i(N, A, B).$   \footnote{ While it is not necessary for the rest of
the paper, here is the algebro-geometric interpretation of the Toda
lattice equations above.  The curve $\CalB_{A,B}$ is a double cover of
the $\lambda$ line. It can be compactified by adding two points $P$ and
$Q$ over $\lambda=\infty$  and $\lambda^{-1}=z$ is a parameter at these
two points.  Further, the function $\alpha$ has a pole of order $N$ at
$P$ and a zero of order $N$ at $Q.$  Let $h(z)=\alpha z^N,$ regarded as
a function of $z$ in a neighborhood of $P$.   To simplify notation,
let's assume that $h(0)=1.$  If $A$ and $B$ are generic with $h(0)=1,$
then given a generic effective  divisor $D$ of degree $g,$  the genus
of $\CalB_{A,B}$, we can find unique functions $s_{n}$ in $L(D+n(P-Q))$
 so that the top coefficient of the Laurent expansion of $s_{n}$ is
one. Thus $s_{0}=1$ and $s(N)=\alpha.$ We define $d_{i}(n)$ to be the
$i^{th}$ Taylor coefficient of $z^n s_{n}$ at $P.$ Thus  $d_i(N)$  is
the $i^{th}$ Taylor coefficient of $h(z)$ expanded as a power
series in $z.$ In particular, $d_i(N)$ does not depend on $D,$ but only
on the curve $\CalB_{A,B}$. It is an easy consequence of Riemann-Roch
that $s_{n}$, $s_{n+1}$, $s_{n-1}$, and  $\lambda s_{n}$ satisfy a
linear relation, which turns out to be of the form  $$\lambda
s_{n}+s_{n+1}=A_{D}s_{n}+B_{D}s_{n-1}$$    where $A_{D}$ and $B_{D}$
are functions which depend only on the linear equivalence class of $D.$
Further, there is a divisor $D_{0}$ so that  $A_{D_{0}}$ is the
original $A$ and $B_{D_{0}}$ is the original $B.$  It is easy to find a
recursive formula for the $d_{i}(n)$ in terms of the $A_{D}$ and
$B_{D}$. We can map the curve $\CalB_{A,B}$ to the Jacobian $J_g$ of
$\CalB_{A,B}$ by $j(x)= D_{0}+P-x.$ The derivative of $j$ at $P$ is an
element of the tangent space of $J_g$ at $D_0,$ which can be extended
to a translation invariant vector field on $J_g.$ The first Toda
equation is just the derivatives of the functions $A_{D}$ and $B_{D}$
as functions on $J_g$ with respect to this translation invariant vector
field.  The higher Toda equations come in a similar way.  Instead of
using the tangent vector to $j(\CalB_{A,B}),$ one uses vectors in the
higher osculating subspaces of $j(\CalB_{A,B})$ at $P.$ The calculation
of the Toda equations from this point of view is essentially worked out
in \cite{1,\S 3}.  The Toda equations are the case $M=1$ of that paper.
 The paper shows how to calculate the third Toda equation, but the
other Toda equations are similarly worked out.}

Given $g\in \CalC,$ our question is whether it is possible to find $N$ and
$(A,B) \in \CalC_N$ so that the Toda flows with initial data $A$
and $B$   look like the KdV flows for $g$ and so that $\CalB_g$  looks
like $\CalB_{A,B}.$  Rather than define precisely what we mean by {\it
looks like}, we will just try to solve the problem through a series of
attempts. We first define $$\epsilon={1\over N}.$$ We will use this
notation consistently throughout the paper. We first consider the
problem of making $\CalB_{A,B}$ look like $\CalB_g.$  We try the
following: $$A(n) = 2+\epsilon^2 g(\epsilon n)$$ and $$B(n)=-1.$$
Suppose that $L_g(\psi)=\lambda \psi$ and $T(\psi)=\alpha \psi$ so that
$(\lambda,\alpha)\in \CalB_g.$ Define $\phi(n)=\psi(n \epsilon).$
Plugging this $A$ and $B$ into
$$L_{(A,B)}(\phi)(n)={-\phi(n+1)+A(n)\phi(n)+B(n)\phi(n-1)\over \epsilon^2}$$
and using
Taylor's theorem to expand $\psi(\epsilon (n+k))$ around $\epsilon n,$
we see that $$L_{A,B}(\phi)(n)= -{d^2\over dx^2}\psi(\epsilon
n)+g(\epsilon n) \psi(\epsilon n)+O(\epsilon ).  $$  Since
$L_g(\psi)=\lambda \psi,$ we have $$L_{A,B}(\phi)(n)= \lambda
\phi(n)+O(\epsilon ).  $$  and $T(\phi)=\alpha \phi.$  Thus $\phi$ is
an approximate eigenfunction of $L_{A,B}$ and $T.$  It is at least
plausible that there is a true eigenfunction $\phi_\epsilon$ of
$L_{A,B}$ and $T$ with  $L_{A,B}(\phi_\epsilon)=\lambda_\epsilon
\phi_\epsilon$ and $T(\phi_\epsilon)=\alpha_\epsilon \phi_\epsilon$ and
$$\lim_{N\to \infty} (\lambda_\epsilon,\alpha_\epsilon
)=(\lambda,\alpha).$$ Conversely, if $\lambda$ is fixed  and $N$ is
large, it is
reasonable to assume that  if we have $(\lambda,\alpha)\in \CalB_N$,
then there are points of $\CalB_g$ nearby.        Thus thus $\CalB_g$
looks like  $\CalB_{A,B}$ with this  choice of $(A,B).$  More precisely,
we conjecture that
if $K$ is any constant and $$S_K=\{(\lambda,\alpha): |\lambda|<K\},$$
we can make the two sets $S_K\cap \CalB_g$ and $S_K\cap \CalB_{A,B}$ be
as close together as we want by making $N$ large.

The  problem with this
approach is
that the above choice of $A$ and $B$ is somewhat arbitrary.  If we let $$A(n) =
2+\epsilon^2 (\frac 1 2 g(\epsilon n)+q(\epsilon n))$$ and
$$B(n)=-1+	 \epsilon^2
(\frac 1 2 g(\epsilon n)-q(\epsilon n)),$$ the same argument shows that
$\CalB_g$
looks like  $\CalB_{A,B}$ with this choice of
$(A,B).$  However, the first conserved quantity
$$d_1(N)=\sum_{k=0}^{N-1} A(k)=2N+\epsilon\int_0^1 \frac 1 2
g(x)+q(x)\,dx+
O(\epsilon^2)$$
depends on our choice of $q.$  That is, we know what $A+B$ is to order
$\epsilon^3$, but we do not know how to define $A$ and $B.$  But the
conserved quantities depend on $A$ and $B,$ not just on $A+B.$
So given $g,$ we want to have a canonical choice of $A$ and $B$ in terms
of $g$.  Then we can hope that the conserved quantities of the $A$ and
$B$ will just depend on $g.$  However, we want this canonical choice of
$A$ and $B$ to be in some sense invariant under the various Toda flows.

We will first concentrate on  the first Toda equation: $$D_{1,1}(A,B)(n)=
(B(n)-B(n+1))N$$ and $$D_{2,1}(A,B)(n)=B(n)(A(n)-A(n-1)N.$$
Let us  consider two functions $f_t(x)=f(x,t)$ and $h_t(x)=h(x,t)$.  Set
$$A_t(n)=2+\epsilon^2(f(\epsilon n,t)+h(\epsilon n,t))$$ and
$$B_t(n)=-1+\epsilon^2(f(\epsilon n,t)-h(\epsilon n,t))$$ and assume
that $f(x,0)=g(x)/2.$ With this choice of $A$ and $B$, $\CalB_g$ still looks
like the Bloch spectrum attached to $A_0$ and $B_0.$ In general,
$\CalB_{2f_t}$ looks like $\CalB_{A_t,B_t}.$ Let's consider the
equations  $${d\over dt}A_t(n)=D_{1,1}(A,B)(n)+O(\epsilon^3)$$ and
$${d\over dt}B_t(n)=D_{2,1}(A,B)(n)+O(\epsilon^3).$$ That is we ask that
$A_t(n)$ and $B_t(n)$ be approximate solutions to the Toda flow
problem.  This means that  $\CalB_{A_t,B_t}$ presumably do not move
much as $t$ changes and hence $\CalB_{2f_t}$ also does not change much
with $t.$ Plugging in and using Taylor's theorem,   we get $${d\over
dt}(f+h)(\epsilon n)=-{d\over dx}(f-h)(\epsilon n)+O(\epsilon)$$ and
$${d\over dt}(f-h)(\epsilon n)=-{d\over dx}(f+h)(\epsilon
n)+O(\epsilon).$$  So if $f$ and $h$ are going to produce solutions to
the Toda flow that are even approximately correct, then we should ask
that  $${d\over dt}(f+h)(x,t)=-{d\over dx}(f-h)(x,t)+O(\epsilon)$$ and
$${d\over dt}(f-h)(x,t)=-{d\over dx}(f+h)(x,t)+O(\epsilon).$$
Obviously, we cannot achieve this by taking $f=h,$ which corresponds to
the case of $B=-1.$  However, we can take $h=0$ and then we do get an
approximate solution to the Toda flow with the choice
$$ A_t(n)=2+\epsilon^2(f(\epsilon n,t))\tag1$$ and
$$ B_t(n)=-1+\epsilon^2(f(\epsilon n,t))\tag2$$ with $${d\over
dt}f(x,t)=-{d\over dx} f(x,t).$$   Equations (1) and (2) define a
canonical choice for $A$ and $B$ in terms of $f=g/2$ modulo $\epsilon^3.$
Of course, this is a rather uninteresting solution.

(1) and (2) only define $A$ and $B$ to order $\epsilon^3.$  We next ask
if we can modify (1) and (2) to be canonical to order $\epsilon^4.$
$$  A_t(n)=2+\epsilon^2 f(\epsilon n,t)-
\epsilon^3\Phi_{3}(f)(\epsilon n,t)\tag4$$
and
$$  B_t(n)=-1+\epsilon^2 f(\epsilon n,t)+
\epsilon^3\Phi_{3}(f)(\epsilon n,t),\tag5$$
where $\Phi_{3}(f)$ is an unknown polynomial in $f$ and the derivatives of
$f.$  We want to choose $\Phi_{3}(f)$ in such a way that
$${d \over dt}A_t=
(B(n)-B(n+1))N+O(\epsilon^4)$$ and
$${d \over dt}B_t=B(n)(A(n)-A(n-1))N+O(\epsilon^4)$$ for some choices at
least of $f(x,t).$  It is not obvious that such a $\Phi_{3}$ should
exist.  Using Taylor's theorem to evaluate $A(n+i)$ and
$B(n+i),$ we should have
$$ {df \over dt}-{d\Phi_{3}(f)\over dt}\epsilon = -{df \over
dx}+\epsilon(-{d^2f\over 2dx^2}-{d\Phi_{3}(f)\over dx})\tag6 $$
and
$$ {df \over dt}+{d\Phi_{3}(f)\over dt}\epsilon = -{df \over dx}
+\epsilon({d^2f\over 2dx^2}+{d\Phi_{3}(f)\over dx}).\tag7$$
 In particular, adding (6) and (7), we see that $f(x,t)$ should
satisfy
$${d\over  dt}f(x,t)=-{d\over dx} f(x,t)+O(\epsilon^2).\tag8$$
Now since $${d\over  dt}f(x,t)=-{d\over dx} f(x,t)+O(\epsilon),$$
 we have $${d\Phi_3(f)\over dt} = -{d\Phi_3(f)\over dx}+O(\epsilon)$$ by
using the chain rule to evaluate $${d\Phi_3(f)\over dt}.$$
So assuming that (8) holds, then (6) will hold if $${d\Phi_3(f)\over
dx}=-\frac 1 4 {d^2f\over dx^{2}}.$$  Thus we can take  $$\Phi_3(f)
=-\frac 1 4 {df\over dx}.$$  Thus if (8) holds and $$
A_t(n)=2+\epsilon^2(f(\epsilon n,t))-
\epsilon^3(-\frac 1 4 {df\over dx}(\epsilon n,t))$$
and
$$  B_t(n)=-1+\epsilon^2(f(\epsilon n,t))+
\epsilon^3(-\frac 1 4 {df\over dx}(\epsilon n,t)),$$ then $A$ and $B$
satisfy the Toda equations to order $\epsilon^4.$
Next we attempt to continue this process.   That is, we try to find a
function $$\Phi_4(f)$$ of $f$ and the derivatives of $f$ so that if we
define  $$
A_t(n)=2+\epsilon^2(f(\epsilon n,t))-
\epsilon^3(-\frac 1 4 {df\over dx}(\epsilon n,t))-
\epsilon^4\Phi_4(f)(\epsilon n,t)\tag9$$
and
$$  B_t(n)=-1+\epsilon^2(f(\epsilon n,t))+
\epsilon^3(-\frac 1 4 {df\over dx}(\epsilon n,t))
+\epsilon^4\Phi_4(f)(\epsilon n,t),\tag10$$
then $A$ and $B$ satisfy the Toda equations to order $\epsilon^5,$
provided that $f(x,t)$ satisfies some condition.    To find this
condition on $f,$ we plug (9) and (10) into
$${d \over dt}A_t=
(B(n)-B(n+1))N+O(\epsilon^5)\tag11$$ and
$${d \over dt}B_t=B(n)(A(n)-A(n-1))N+O(\epsilon^5).\tag12$$ Adding (11)
and (12), we see that $f$ must satisfy:
$${d\over  dt}f(x)=-{d\over dx} f(x)+(-\frac 1 {24} \dthree + \frac 1
2 f(x)
\done )\epsilon^2+O(\epsilon^3).\tag13$$
Assuming (13) and using the fact that
$${d\Phi_4(f)\over dt} = -{d\Phi_4(f)\over dx}+O(\epsilon),$$ some
computation shows that
$${d\Phi_4(f)\over dx}=-\frac 1 4 f{df\over dx}.$$  So we see that if
$${\Phi_4(f)}=-\frac 1 8 f^2,$$  then both (11) and (12) will hold if
(13) holds.

When we say that $f$ satisfies (13), we mean that (13) is satisfied in
the asymptotic sense.  More generally, consider the equation
$${df \over dt}=c {df \over dx} +\epsilon\Psi_{1}(f)+\ldots
+\epsilon^n\Psi_{n}(f)+O(\epsilon^{n+1}),\tag14$$
where the $\Psi_{k}$ are polynomials in $f$ and the $x$-derivatives of $f.$
Of course, it can be extremely difficult to solve $${df \over
dt}=c{df\over dx}
+\epsilon\Psi_{1}(f)+\ldots
+\epsilon^n\Psi_{n}(f).\tag15$$ So instead we try an asymptotic solution
of the form
$$f(x,t)=f_{0}(x,t) +f_{1}(x,t)\epsilon+\ldots +f_{n}(x,t)\epsilon^n.\tag16$$
To solve our problem, substitute (15) into (16) and multiply everything
out and denote  the coefficient of $\epsilon^k$ by $V_{k}$.  Then we can
easily find $f_{k}$ for $k$ from 0 to $n$ so that $V_{k}=0$
for $k$ from 0 to $n.$ Indeed, all we have to do is to solve
$${df_{k}\over dt}=c{df_{k}\over dx}+L_{k}(f_0,\ldots,f_{k-1}),$$ where the
$L_{k}$ are polynomials in the known $f_0,\ldots,f_{k-1}$ and their
derivatives.  Given initial
conditions, we can then easily solve these equations. So an asymptotic
solution in this sense to (13) will solve (11) and (12) to order
  $\epsilon^5$ for $t$ in any given bounded region.  That is, for any
given $T$, there is an integer $N$ and a constant $C$ so that if $n>N,$
then $$\vert {d \over dt}A_t-
(B(n)-B(n+1))N\vert < C\epsilon^5 $$ and
$$\vert{d \over dt}B_t-B(n)(A(n)-A(n-1))N\vert<C\epsilon^5$$ if $\vert t
\vert < T.$  Of
course, these estimates may break down over long time intervals, i.e.
$C$ may depend on $T.$  However, this kind of asymptotic analysis
provides a tool for guessing the behavior of the equations of the Toda
hierarchy, although rigorously connecting the asymptotic analysis and
the behavior of the equations is probably extremely difficult.

I conjecture that this process of finding $\Phi_k$ can be continued
indefinitely.  At any
rate suppose we define $$\eqalign{&R(f)=-\frac 1 4 \done-\frac 1 8 f(x)^2
\epsilon +\frac 1 {192}
 \dthree \epsilon^2+\cr &\qquad
 (\frac 1 {64} f(x) \dtwo+
\frac 1 {64} \done^2-1/32 f(x)^3)\epsilon^3\cr & \qquad+
 (-\frac 1 {7680} \dfive+\frac 1 {64}f(x)^2
\done)\epsilon^4\cr &\qquad+
 (\frac 3 {256} \done^2 f(x)+\frac 3 {512} f(x)^2 \dtwo-\frac 5 {512} f(x)^4
\cr &\qquad-
\frac 1{1536} f(x) \dfour-\frac 3{2048} \dtwo^2-\frac 1{384} \done
\dthree)
\epsilon^5}$$
and $f$ satisfies $$\eqalign{&{df\over dt}=
-\done+(-\frac 1 {24} \dthree + \frac 1 2 f(x)
\done )\epsilon^2\ski
(-\frac 1 {1920} \dfive + \frac 1 {16} f(x)^2 \done +\frac 1 {32}
\done \dtwo+
 \frac 1 {48} f(x) \dthree )\epsilon^4
\ski  ( \frac 1 {128} \dtwo^2 +\frac 1 {128} \done \dthree
)\epsilon^5
\ski
(-\frac 1 {322560} \dseven + \frac 1 {128} \done ^3-\frac 1 {1536} \dtwo
\dthree \cr &\qquad+\frac 1 {384} f(x)^2 \dthree+
\frac 1 {64} \done f(x) \dtwo +\frac 1 {64} f(x)^3 \done \cr &\qquad+
\frac 1 {3840}
f(x) \dfive + \frac 1 {1536} \done \dfour)\epsilon^6
}$$
If we define    $$  A_t(n)=2+\epsilon^2(f(\epsilon n,t))-
\epsilon^3R(f)(\epsilon n,t)\tag17$$
and
$$  B_t(n)=-1+\epsilon^2(f(\epsilon n,t))+
\epsilon^3R(f)(\epsilon n,t),\tag18$$
then $${d \over dt}A_t=
-(B(n+1)-B(n))N+O(\epsilon^9)$$ and
$${d \over dt}B_t=B(n)(A(n)-A(n-1))N+O(\epsilon^9).$$  These computations
were carried out using Maple, as is described in the Appendix.

Next let's consider the second pair of equations in the Toda hierarchy.
Actually, the answer is more comprehensible if we try to solve the equations:
 $$(D_{1,2}+2D_{1,1})(A_t,B_t)={d\over dt}A_t+O(\epsilon^9)\tag19$$ and
$$(D_{2,2}+2D_{2,1})(A_t,B_t)={d\over dt}B_t+O(\epsilon^9).\tag20$$  Now define
$A_{t}$ and $B_{t}$ by the equations (17) and (18). In order for (19)
and (20) to hold, we must have that
$$\epsilon^2{df\over
dt}=(D_{1,2}+2D_{1,1}+D_{2,2}+2D_{2,1})(A_t,B_t)+O(\epsilon^9).$$
Plugging in our definition of $A$ and $B$ and using Taylor's theorem
many times, we see we must have

$$\eqalign{&{df\over dt}= (-\frac 1 {4} \dthree + 3 f(x)
\done )\epsilon^2\cr &\qquad+
(\frac 3 8 f(x) \dthree - \frac 9 8 f(x)^2 \done - \frac 1 {64} \dfive
\ski
\frac {11}{16} \done \dtwo )\epsilon^4\ski
  ( \frac 3 {64} \dtwo^2 +\frac 3 {64} \done \dthree
)\epsilon^5
\ski
((-\frac 1 {2560} \dseven - \frac 1 {64} \done ^3+\frac {49} {768}
\dtwo
\dthree \cr &\qquad-\frac {5} {64} f(x)^2 \dthree-
\frac 7 {32} \done f(x) \dtwo -\frac 9 {32} f(x)^3 \done \cr
&\qquad+\frac {11} {640}
f(x) \dfive + \frac {35} {768} \done \dfour)\epsilon^6
}.\tag21$$
Conversely, if (21) is valid and we use (17) and (18) to define $A$ and
$B,$ then (19) and (20) are valid.

Next we repeat this process for the third set of equations in the Toda
hierarchy.  Here the result is that if
$$\eqalign{&{df\over dt}=(\frac 5 4 f(x) \dthree - \frac {15} 2 f(x)^2 \done -
\frac
1 {16} \dfive
\ski
\frac {5}{2} \done \dtwo )\epsilon^4\ski
((-\frac 1 {192} \dseven - \frac 5 {4} \done ^3+\frac {155} {192}
\dtwo
\dthree \cr &\qquad-\frac {45} {32} f(x)^2 \dthree-
\frac {85} {16} \done f(x) \dtwo +\frac {15} {8} f(x)^3 \done \cr
&\qquad+\frac {11} {64}
f(x) \dfive + \frac {47} {96} \done \dfour)\epsilon^6
}$$
and we use (17) and (18) to define $A$ and
$B,$ then
$$(D_{1,3}-2D_{1,1}+2D_{1,2})(A_t,B_t)={d\over dt}A_t+O(\epsilon^9)\tag19$$ and
$$(D_{2,3}-2D_{2,1}+2D_{2,2})(A_t,B_t)={d\over dt}B_t+O(\epsilon^9)\tag20.$$

Finally, for the fourth Toda equation, we have that
$$(D_{1,4}+4D_{1,1}-2D_{1,2}+2D_{1,3})(A_t,B_t)={d\over
dt}A_t+O(\epsilon^9)\tag19$$ and
$$(D_{2,4}+4D_{2,1}-2D_{2,2}+2D_{2,3}))(A_t,B_t)={d\over
dt}B_t+O(\epsilon^9)\tag20$$
if $$\eqalign{&{df\over dt}=
((-\frac 1 {64} \dseven - \frac {35} {8} \done ^3+\frac {35} {16}
\dtwo
\dthree \cr &\qquad-\frac {35} {8} f(x)^2 \dthree-
\frac {35} {2} \done f(x) \dtwo +\frac {35} {2} f(x)^3 \done \cr
&\qquad+\frac {7} {16}
f(x) \dfive + \frac {21} {16} \done \dfour)\epsilon^6
}$$

This suggests the following conjecture:

\proclaim{Conjecture}
  There are polynomials $\Phi_k(f)$ and $\Psi_{k,j}(f)$
in $f$ and the derivatives of $f$ so that if
$$R(f)=\sum_{k=1}^{L}\Phi_{k}(f)\epsilon^{k}$$ and
$$Z_{j}(f)=\sum_{k=0}^{L+1}\Psi_{k,j}(f)\epsilon^k
$$and we define
$$  A_t(n)=2+\epsilon^2(f(\epsilon n,t))-
\epsilon^3R(f)(\epsilon n,t)$$
and
$$  B_t(n)=2+\epsilon^2(f(\epsilon n,t))+
\epsilon^3R(f)(\epsilon n,t)$$
and $f$ satisfies the equation
$${df\over dt}= Z_{j}(f)+O(\epsilon^{L+1}
),$$ then $A_t(n)$ and $  B_t(n)$ satisfy the $j^{th}$ Toda equations to
order $\epsilon^{L+4}.$ Further, by taking suitable linear combinations of
the $Z_{j}(f),$  we can produce asymptotic series whose leading terms in
$\epsilon $
are the KdV hierarchy if $L$ is large enough.
\endproclaim

For the rest of the paper, let us assume this conjecture is true.
Then we can also calculate out the conserved quantities asymptotically.
First, let's consider $d_1(N),$ the first conserved quantity.  We have
$$d_1(N)=\sum_{n=0}^{N-1}A(n).$$  By the Euler-Maclaurin summation
formula, we have that asymptotically
$$d_1(N)= N (\int_{0}^{1}(2+\epsilon^2 f(x)-\epsilon^3
R(f))\,dx)
.$$
So we can write out the first few terms of the expansion of  $d_1(N):$
$$d_{1}(N)=\frac 1 \epsilon  (2+ \epsilon \int_{0}^{1}f(x)\,dx +\frac
1 8 \epsilon^3 \int_{0}^{1}f(x)^2\,dx+\frac 1 {32}  \epsilon^5
\int_{0}^{1}f(x)^3\,dx +\ldots. $$
In particular, the following quantity is conserved asymptotically:
$$C1=\int_{0}^{1}f(x)\,dx +\frac
1 8 \epsilon^2 \int_{0}^{1}f(x)^2\,dx+\frac 1 {32}  \epsilon^4
\int_{0}^{1}f(x)^3\,dx +\ldots. $$
Next consider $d_2(N).$ We can write $$d_2(N)=
(\sum_{n=0}^{N-1}A(n))^2+\sum_{n=0}^{N-1}B(n)+\sum_{n=0}^{N-1}A(n)^2.$$
Again the Euler-Maclaurin summation formula allows us to calculate each
of the terms as an integral.  In order to make the formulas look nicer,
we consider to conserved quantity
$$C2 =-\frac 4 3 (d_2(N)-\frac 2 {\epsilon^2} + \frac 3
{\epsilon^2})+\frac {2-\epsilon} \epsilon C1) - \frac 1 2 C1^2.\tag 22 $$
Then we have an expansion for $C2$
$$\epsilon^3\int_0^1 f(x)^2\,dx+\epsilon^5(\frac 1 4 \int_0^1 f(x)^3\,dx
-\frac 1 {24} \int_0^1 f(x) {d^2 f(x) \over dx^2} \,dx)+o(\epsilon^7) .$$
Thus if we solve a given Toda equation asymptotically to high enough
order, we can expect $$\int_0^1 f(x)^2\,dx+\epsilon^2(\frac 1
4 \int_0^1 f(x)^3\,dx
-\frac 1 {24} \int_0^1 f(x) {d^2 f(x) \over dx^2} \,dx)+\ldots  $$ to be
conserved to order $\epsilon^4.$

The third conserved quantity is $C3$ which has the expression:
$$d_3(N)-\frac 1 {12\epsilon}(-16+72\epsilon-56\epsilon^2-24 C1 \epsilon+60 C1
\epsilon^2$$
$$-36 C1 \epsilon^3 + 12 C1^2 \epsilon^2+18 C2
\epsilon^2+2C1^3\epsilon^3+18C2\epsilon^3-9 C1 C2 \epsilon^3). $$
We have the following expression for $C3$
$$ \epsilon^5(-\frac 7 {12} \int_0^1 f(x)^3\,dx
+\frac 1 {8} \int_0^1 f(x) {d^2 f(x) \over dx^2} \,dx)+\ldots $$

Let's use (19) and (20) to try to find a numerical method for solving KdV.
We first define $a(n)=N^2(A(n)-2)$ and $b(n)=N^2(B(n)+1).$ Then we define
$$L_{a,b}(k)=2 b(k + 1) - 2 b(k) - a(k + 1) + a(k - 1)$$ and
$$M_{a,b}(k)=2 a(k) - 2a(k - 1) - b(k + 1) + b(k - 1).$$
 $$F_{a,b}(k)= b(k+1)a(k)+b(k+1)a(k+1)-b(k)a(k)-b(k)a(k-1)$$ and
$$ \displaylines{G_{a,b}(k)=-2b(k)a(k)+2b(k)a(k-1)+a(k)^2-a(k-1)^2+
b(k)b(k+1)-b(k)b(k-1)\cr+
(-b(k)a(k)^2+b(k)a(k-1)^2)/N^2.}$$

(19) and (20) become
$${da\over dt}(k)=N(  L_{a,b}(k)+ \epsilon^2 F_{a,b}(k))\tag 23$$
and $${db\over dt}(k)=N(  M_{a,b}(k)+ \epsilon^2 G_{a,b}(k)).\tag 24$$
Let's suppose we have a function $f(x)$ we want to use as initial data for
KdV where $f$ is periodic with period 1. Define $c(n)=f(n/N).$  Then we can
mimic  (17) and (18) to define $a(n)$ and $b(n)$ to any desired degree
of accuracy in $\epsilon.$  Thus let's define $$\Delta (b) (k) ={
N(b(k+1)-b(k-1))\over 2}.$$ So we can approximate the first derivative
of $f$ by $$\Delta(c) -\frac 1 6 \epsilon^2 \Delta^{3}(c).$$

For instance, we could set
$$a(n)= c(n)+\frac 1  4 \epsilon^2 (\Delta (c)(n) -\frac 1 6 \epsilon^2
\Delta^{3}(c))
+\frac 1 8{\epsilon^3}
c(n)^2 -\frac 1 {192}\epsilon^4\Delta^3(c).\tag 25$$ and
$$b(n)= c(n)-\frac 1  4 \epsilon^2 (\Delta (c)(n) -\frac 1 6 \epsilon^2
\Delta^{3}(c))
-{\epsilon^3\over 8}
c(n)^2 +\frac 1 {192}\epsilon^4\Delta^3(c))\tag 26$$ and use these $a$ and $b$
as initial conditions   for (23) and (24) to get functions $a(n,t)$
and $b(n,t).$  We can reasonable expect that $$a(n,t)+b(n,t)\over 2$$ will
remain close to $f(n/N,t),$  where $f(x,t)$ is a solution to $$
{df\over dt}= (-\frac 1 {4} \dthree + 3 f(x)
\done )\epsilon^2. $$  In fact, this method should be accurate to order
$\epsilon^4.$  Further, we can expect a whole series of conserved
quantities  similar to (22).

There are serious problems in discretizing the equations (23) and (24)
in time efficiently.  This problem can be seen by setting
$a(n)=\epsilon^2 v(n/N)$ and $b(n)=\epsilon^2 w(n/N).$ Then to order
$\epsilon,$  we have  $${dv\over dx} = 2{dw\over dx} -2{dv\over dx}
\tag 27$$ and $${dw\over dx} = 2{dv\over dx} -2{dw\over dx}.\tag 28$$
So while arbitrary solutions tend to move with speeds on the order of
magnitude 1, the special solutions we have constructed with (25) and
(26) tend to move with speed $\epsilon^2.$  Thus if one applies some
explicit method such as Runge-Kutte, one is virtually guaranteed that
the differencing in time will produce instability if the time step size
is not small with respect to 1.  One can intuitively understand  the
situation as follows:  the equations (23) and (24) define a vector
field on the space $V$ of functions $a$ and $b.$  For most of  the $a$
and $b$ coming from smooth  data, the vectors in this field have
order of magnitude one.  But for the set $W$ of those coming from (25)
and (26), the order of the magnitude of these vectors is $\epsilon^2.$
Thus $W$ is almost fixed by the flow. In fact, $W$ tends to act as a
neutral fixed set, so solutions nears $W$ tend to stay near $W.$ A large step
size will tend to
push away from $W,$ and as more steps  are taken, the distance from $W$
increases exponentially.  In practice, the high frequencies are most
increased and so instability results.  I believe that this problem can
be overcome using the Crank-Nicholson method.  Given $a$ and $b$ we wish
to find $a''$ and $b'',$  the results of propagating $a$ and $b$ using
(23) and (24) over a time step $\Delta t.$  We try to solve the
equations  $$a(k) + {\Delta t\over 2}(N^3 L_{a,b}(k)+ N
F_{a,b}(k))=a'(k)-{\Delta t\over 2}(N^3
L_{a',b'}(k)+ N F_{a',b'}(k))$$ and
  $$b(k) +{\Delta t\over 2}(N^3 M_{a,b}(k)+ N G_{a,b}(k))=b'(k)-{\Delta
t\over 2}(N^3 M_{a',b'}(k)+ N G_{a',b'}(k)).$$   to obtain an
Crank-Nicholson approximations $a'$ and $b'$ to $a''$ and $b''.$
If we take the step size small with respect of $N^2,$  then $${\Delta
t\over 2} N F_{a,b}(k)$$ is small.  So in solving for $a'$ and $b'$ in
terms of $a$ and $b,$  the non-linear terms are small and so one should
be able to find an algorithm to solve these equations using some version
of Newton's method.  At any rate, I tried a rather {\it ad hoc} method
which was stable on the examples I tried.

  \head Appendix \endhead

The calculations in this paper were made using the computer algebra
system Maple V, Release 2. Here is an explanation of the functions in
this program, followed by the program:

$\bullet$ dt is a procedure designed to differentiate $g$ by $t,$ given
that $${df\over dt}=h.$$ So dt acts as a derivation of sums, products
and exponentials.  Further, $${d \over dt}{d\over dx} $$ is replaced by
$${d \over dx}{d\over dt}. $$\medskip
$\bullet$ d(i,n) is $d_i(n).$  It implements the recursive relation
defining the $d_i(n).$
\medskip
$\bullet$ translate produces the Taylor expansion of $f(x+n\epsilon)$ to
order $N$.
\medskip
$\bullet$ R is defined in the main part of this paper.

\medskip
$\bullet$ aa implements the recursive definition of the Toda equations.
\medskip
$\bullet$
XZ[p] is $D_{1,k}$ and YZ[p] is $D_{2,k}.$
\medskip
$\bullet$
AA[p] and XX are the derivative of $f$ for that particular Toda equation.
\medskip
$\bullet$
check[p] checks that the time derivative of $B$ is the same as the
derivative calculated using the $p^{th}$ Toda equation.
\medskip
$\bullet$ The last three lines calculate the various linear combinations
of AA[p] to get the KdV hierarchy.
\vskip .4in
\NoBlackBoxes
{\eighttt
\  \  \  \  dt\  :=\  proc(g,f,h)

 \  \  \  \  \  \  \  \  \  local\  i,n,p;

 \  \  \  \  \  \  \  \  \  \  \  \  \  if\  type(g,`+`)\  then

 \  \  \  \  \  \  \  \  \  \  \  \  \  \  \  \  \  p\  :=\  0;

 \  \  \  \  \  \  \  \  \  \  \  \  \  \  \  \  \  n\  :=\  nops(g);

 \  \  \  \  \  \  \  \  \  \  \  \  \  \  \  \  \  for\  i\  to\  n\  do\  \
p\  :=\  p+dt(op(i,g),f,h)\  od;

 \  \  \  \  \  \  \  \  \  \  \  \  \  \  \  \  \  RETURN(p)

 \  \  \  \  \  \  \  \  \  \  \  \  \  elif\  type(g,`*`)\  then

 \  \  \  \  \  \  \  \  \  \  \  \  \  \  \  \  \  p\  :=\  0;

 \  \  \  \  \  \  \  \  \  \  \  \  \  \  \  \  \  for\  i\  to\  nops(g)\
do\  \  p\  :=\  p+dt(op(i,g),f,h)/op(i,g)\  od;

 \  \  \  \  \  \  \  \  \  \  \  \  \  \  \  \  \  RETURN(expand(g*p))

 \  \  \  \  \  \  \  \  \  \  \  \  \  elif\  type(g,`**`)\  then

 \  \  \  \  \  \  \  \  \  \  \  \  \  \  \  \  \
RETURN(op(2,g)*op(1,g)**(op(2,g)-1)*dt(op(1,g),f,h))

 \  \  \  \  \  \  \  \  \  \  \  \  \  elif\  type(g,function)\  and\
op(0,g)\  =\  f\  then\  RETURN(h)

 \  \  \  \  \  \  \  \  \  \  \  \  \  elif\  type(g,function)\  and\
op(0,g)\  =\  diff\  then

 \  \  \  \  \  \  \  \  \  \  \  \  \  \  \  \  \
RETURN(diff(dt(op(1,g),f,h),x))

 \  \  \  \  \  \  \  \  \  \  \  \  \  elif\  type(g,numeric)\  or\  g\  =\
epi\  then\  RETURN(0)

 \  \  \  \  \  \  \  \  \  \  \  \  \  fi

 \  \  \  \  \  \  \  \  \  end\  \  \  \  \  \  \  \  \  \  \  \  \  \  \  \
\  \  \  \  \  \  \  \  \  \  \  \  \  \  \  \  \  \  \  \  \  \  \  \  \  \  \
 \  \  \  \  \  \  \  \  \  \  \  \  \  \  \  \  \  ;

 \  \  \  d\  :=

 \  \  \  \  \  \  proc(i,n)

 \  \  \  \  \  \  options\  remember;

 \  \  \  \  \  \  \  \  \  \  if\  i\  <\  0\  then\  RETURN(0)

 \  \  \  \  \  \  \  \  \  \  elif\  i\  =\  0\  then\  RETURN(1)

 \  \  \  \  \  \  \  \  \  \  elif\  n\  =\  0\  then\  RETURN(0)

 \  \  \  \  \  \  \  \  \  \  elif\  n\  <\  0\  then\
RETURN(d(i,n+1)-A(n)*d(i-1,n)-B(n)*d(i-2,n-1))

 \  \  \  \  \  \  \  \  \  \  elif\  0\  <\  n\  then\
RETURN(d(i,n-1)+A(n-1)*d(i-1,n-1)+B(n-1)*d(i-2,n-2))

 \  \  \  \  \  \  \  \  \  \  fi

 \  \  \  \  \  \  end

 \  \  \  \  \  \  ;

 \  \  \  translate\  :=

 \  \  \  \  \  \  \  proc(f,n)

 \  \  \  \  \  \  \  local\  i,y;

 \  \  \  \  \  \  \  \  \  \  \  y\  :=\  f(x);\  for\  i\  to\  N\  do\  \
y\  :=\  y+diff(f(x),x\  \$\  i)*(epi*n)**i/i!\  od

 \  \  \  \  \  \  \  end\  \  \  \  \  \  \  \  \  \  \  \  \  \  \  \  \  \
\  \  \  \  \  \  \  \  \  \  \  \  \  \  \  \  \  \  \  \  \  \  \  \  \  \  \
 \  \  \  \  \  \  \  \  \  \  \  \  \  \  \  \  \  \  \  \  \  ;

 \  \  \  R\  :=

 \  \  \  proc(x)

 \  \  \  \  \  \  \  -1/4*diff(f(x),x)-1/8*epi*f(x)**2+1/192*diff(f(x),x\  \$\
 3)*epi**2+

 \  \  \  \  \  \  \  \  \  \  \  epi**3*(1/64*f(x)*diff(f(x),x\  \$\
2)+1/64*diff(f(x),x)**2-1/32*f(x)**3)+

 \  \  \  \  \  \  \  \  \  \  \  epi**4*(-1/7680*diff(f(x),x\  \$\
5)+1/64*f(x)**2*diff(f(x),x))+epi**5*(

 \  \  \  \  \  \  \  \  \  \  \
3/256*diff(f(x),x)**2*f(x)+3/512*f(x)**2*diff(f(x),x\  \$\  2)-5/512*

 \  \  \  \  \  \  \  \  \  \  \  f(x)**4-1/1536*f(x)*diff(f(x),x\  \$\
4)-3/2048*diff(f(x),x\  \$\  2)**2-1/384

 \  \  \  \  \  \  \  \  \  \  \  *diff(f(x),x)*diff(f(x),x\  \$\  3))

 \  \  \  end

 \  \  \  ;

 \  \  \  N\  :=\  11;

 \  \  \  F\  :=\  proc(x)\  2+epi**2*f(x)-epi**3*R(x)\  end;

 \  \  \  G\  :=\  proc(x)\  -1+epi**2*f(x)+epi**3*R(x)\  end;

 \  \  \  XZ[1]\  :=\  B(0)-B(1);

 \  \  \  YZ[1]\  :=\  B(0)*(A(0)-A(-1));

 \  \  \  for\  k\  from\  2\  to\  4\  do

 \  \  \  \  \  \  \  aa\  :=\  proc(p,n)

 \  \  \  \  \  \  \  \  \  \  \  \  \  local\  r,X;

 \  \  \  \  \  \  \  \  \  \  \  \  \  \  \  \  \  X\  :=\
d(k-p,n+k)-d(k-p,n);

 \  \  \  \  \  \  \  \  \  \  \  \  \  \  \  \  \  for\  r\  from\  p+1\  to\
k-1\  do\  \  X\  :=\  X-aa(r,n)*d(r-p,n+r)\  od;

 \  \  \  \  \  \  \  \  \  \  \  \  \  \  \  \  \  X

 \  \  \  \  \  \  \  \  \  \  \  \  \  end\  \  \  \  \  \  \  \  \  \  \  \
\  \  \  \  \  \  \  \  \  \  \  \  \  \  \  \  \  \  \  \  \  \  \  \  \  \  \
 \  \  \  \  \  \  \  \  \  \  \  \  \  \  \  \  \  ;

 \  \  \  \  \  \  \  XZ[k]\  :=\  expand(aa(1,0)*B(1)-aa(1,-1)*B(0));

 \  \  \  \  \  \  \  YZ[k]\  :=\  expand(B(0)*(aa(0,0)-aa(0,-1)))

 \  \  \  od;

 \  \  \  A\  :=\  proc(n)\  options\  remember;\  translate(F,n)\  end;

 \  \  \  B\  :=\  proc(n)\  options\  remember;\  translate(G,n)\  end;

 \  \  \  for\  p\  to\  4\  do

 \  \  \  \  \  \  \  XX\  :=\
convert(taylor(1/2*(XZ[p]+YZ[p])/epi**3,epi,12),polynom);

 \  \  \  \  \  \  \  AA[p]\  :=\  collect(expand(XX),epi);

 \  \  \  \  \  \  \  YY\  :=\  taylor(dt(G(x),f,XX)-YZ[p]/epi,epi,9);

 \  \  \  \  \  \  \  YY\  :=\  collect(expand(YY),epi);

 \  \  \  \  \  \  \  check[p]\  :=\  YY;

 \  \  \  \  \  \  \  print(p,check[p])

 \  \  \  od;

 \  \  \  printlevel\  :=\  1;

 \  \  \  AA[1];

 \  \  \  AA[2]\  :=\  collect(expand(AA[2]+2*AA[1]),epi);

 \  \  \  AA[3]\  :=\  collect(expand(AA[3]-6*AA[1]+2*AA[2]),epi);

 \  \  \  AA[4]\  :=\  collect(expand(AA[4]+20*AA[1]-6*AA[2]+2*AA[3]),epi)

}

\enddocument